\documentclass{amsart}
\usepackage{color}
\usepackage[all]{xypic}

\newtheorem{theorem}{Theorem}[section]
\newtheorem{lemma}{Lemma}
\newtheorem{corollary}{Corollary}

\theoremstyle{conjecture}

\theoremstyle{definition}

\theoremstyle{remark}

\theoremstyle{remarks}
\newtheorem*{remarks}{Remarks}
\theoremstyle{example}

\theoremstyle{examples}

\theoremstyle{problems}

\numberwithin{equation}{section}

\begin{document}

\title{Solving for Root Subgroup Coordinates: The SU(2) Case}

\author{Doug Pickrell}
\email{pickrell@math.arizona.edu}

\begin{abstract} In \cite{Pi1} and
\cite{PP} we showed that a loop in a simply connected compact Lie group $\dot K$
has a unique triangular factorization if and only if the loop has a
unique root subgroup factorization (relative to a choice of a reduced sequence of
simple reflections in the affine Weyl group). In this paper we show
that in the $\dot K=SU(2)$ case, root subgroup coordinates are rational functions (with positive denominators)
of the triangular factorization coordinates. We conjecture that in general they are algebraic
functions.



\end{abstract}

\maketitle

\setcounter{section}{-1}

\section{Introduction}\label{Introduction}

In this paper, unless stated otherwise, we suppose that $\dot K:=SU(2)$ and $\dot G:=SL(2,\mathbb C)$.

A triangular factorization
for $g\in L\dot G$ is a factorization of the form
\begin{equation}\label{trifactorization}g(z)=l(z)\cdot m\cdot a\cdot u(z),\end{equation}
where
\[l(z)=\left(\begin{array}{cc}
l_{11}(z)&l_{12}(z)\\
l_{21}(z)&l_{22}(z)\end{array} \right)\in H^0(\Delta^{*},G),\quad
l(\infty )=\left(\begin{array}{cc}
1&0\\
l_{21}(\infty )&1\end{array} \right),\] $l$ has appropriate
boundary values on $S^1$ (depending on the smoothness properties
of $g$), $m=\left(\begin{array}{cc}
m_0&0\\
0&m_0^{-1}\end{array} \right)$, $m_0\in S^1$,
$a=\left(\begin{array}{cc}
a_0&0\\
0&a_0^{-1}\end{array} \right)$, $a_0>0$,
\[u(z)=\left(\begin{array}{cc}
u_{11}(z)&u_{12}(z)\\
u_{21}(z)&u_{22}(z)\end{array} \right)\in H^0(\Delta ,G),\quad
u(0)=\left(\begin{array}{cc}
1&u_{12}(0)\\
0&1\end{array} \right),\] and $u$ has appropriate boundary values
on $S^1$, where $\Delta$ ($\Delta^*$) is the open unit disk
centered at $z=0$ ($z=\infty$, respectively), and $H^0(U)$ denotes
holomorphic functions in a domain $U\subset \mathbb C$.

The basic
fact is that for $g\in L\dot K$ having a triangular factorization,
there is a second unique root subgroup factorization
\begin{equation}\label{coordinate}g(z)=k_1(\eta)^{*}(z)\left(\begin{matrix} e^{\chi(z)}&0\\
0&e^{-\chi(z)}\end{matrix} \right)k_2(\zeta)(z),\quad \vert
z\vert=1,\end{equation} where
\begin{equation}\label{k1product}k_1(\eta)(z)=\lim_{n\to\infty}\mathbf a(\eta_n)\left(\begin{matrix} 1&-\overline{\eta}_nz^n\\
\eta_nz^{-n}&1\end{matrix} \right)..\mathbf
a(\eta_0)\left(\begin{matrix} 1&-\overline{
\eta}_0\\
\eta_0&1\end{matrix} \right),\end{equation} $\mathbf{\chi}(z)
=\sum\mathbf{\chi}_jz^j$ is a $i\mathbb R$-valued Fourier series
(modulo $2\pi i\mathbb Z$),
\begin{equation}\label{k2product}k_2(\zeta)(z)=\lim_{n\to\infty}\mathbf a(\zeta_n)\left(\begin{matrix} 1&\zeta_nz^{-n}\\
-\overline{\zeta}_nz^n&1\end{matrix} \right)..\mathbf
a(\zeta_1)\left(\begin{matrix} 1&
\zeta_1z^{-1}\\
-\overline{\zeta}_1z&1\end{matrix} \right),\end{equation} $\mathbf
a(\cdot )=(1+\vert\cdot\vert^2)^{-1/2}$, and it is understood that
if $g\in C^{\infty}(S^1,K)$, then the coefficients are rapidly
decreasing, and similarly for other function spaces; conversely a root
subgroup factorization as in (\ref{coordinate}) implies that $g$ has a
triangular factorization (\cite{Pi1}).

It is a relatively simple matter to pass from a root subgroup factorization to a triangular
factorization; we will recall how this is done in section \ref{rsfimpliestf}. The main point of
this paper is to explain how to directly find the root subgroup factors $\eta$ and $\zeta$ in terms
of the triangular factorization. The
proof in \cite{Pi1} for the existence of these factors uses (in part) an inverse function theoretic argument.
The basic observation is that the Taylor series centered at $z=0$ for the meromorphic
functions $l_{21}(z)/l_{11}(z)$ and $u_{21}(z)/u_{22}(z)$ have coefficients which are related in a triangular way
with the variables $\eta$ and $\zeta$, respectively. From this we can see that $\eta$ ($\zeta$)
is a rational function of the coefficients of $l$ ($u$) and their conjugates, with positive denominators. Unfortunately we have
not yet found closed form expressions for $eta$ and $\zeta$.

This reveals a surprising fact about the variables $\eta$ and $\zeta$: their natural domain is the set
of all loops $g$ in the formal completion of the complex loop group $ L\dot G$ having a triangular factorization.
For example this immediately implies that (the components of) $\eta$ and $\zeta$ are well-defined random variables
with respect to the invariant measures considered in \cite{Pi0}.

In finite dimensions there is a Gaussian elimination algorithm for finding
the LDU decomposition for a matrix. There does not exist an algorithm for finding the analogous triangular decomposition
for a matrix valued loop - one has to invert a Toeplitz operator. What we are showing is that triangular factorization
and root subgroup factorization are at roughly the same level of complexity.

\subsection{Higher Rank and Generalizations}

Suppose that $\dot K$ is an arbitrary simply connected compact Lie group with simple Lie algebra $\dot{\mathfrak k}$. Fix
a triangular decomposition of the complexified Lie algebra
$$\dot{\mathfrak g}:=\mathfrak k^{\mathbb C}=\mathfrak n_-\oplus \mathfrak h \oplus \mathfrak n_+$$ which is compatible
with $\mathfrak k$ in the sense that $\mathfrak k\cap\mathfrak h$ is maximal abelian.
In this case there is again a triangular factorization for $g:S^1\to K$, $g(z)=l(z)mau(z)$, if and only if there is a root subgroup factorization,
$g(z)=k_1(\eta)^*(z)exp(\chi(z))k_2(\zeta)(z)$, where in the case $rank(\mathfrak k)>1$, the details of the root subgroup factorization
depend additionally on choices of reduced factorization of the longest element $w_0$ of the Weyl group, and a reduced sequence of simple reflections for the affine Weyl group; see \cite{PP}.

Based on experiments with $SU(3)$ and $SU(4)$, it appears that in this more general context, the components of $\eta$ and $\zeta$
are algebraic functions of the triangular factorization coordinates. But solving for the components is far less straightforward.

The result in this paper does resurrect the hope that there might be a way to similarly solve for the root subgroup coordinates
for homeomorphisms of a circle; see \cite{Pi2} (this is what has long motivated me to search for formulas similar
to those in this paper). The formulas that we find also apply for root subgroup factorization
for loops in $G_0=SU(1,1)$, as in \cite{CP}. It is possible that these formulas might help clarify some of the complications for root subgroup factorization that arise in that noncompact context.

\subsection{Notation}

For a function $f:U\subset \widehat{\Sigma}\to \mathcal L(\mathbb C^N)$, define
$f^*(q)=f(R(q))^*$, where $(\cdot)^*$ is the Hermitian adjoint. If $f\in H^0(\Sigma)$ (i.e. a holomorphic function in
some open neighborhood of $\Sigma$), then $f^*\in H^0(\Sigma^*)$. If $q\in S$, then $f^*(q)=f(q)^*$, the ordinary
complex conjugate of $f(q)$.

\section{From Root Subgroup Factorization to a Triangular
Factorization: The $SU(2)$ Case} \label{rsfimpliestf}

Suppose that $g:S^1\to SU(2)$ has a root subgroup factorization as in (\ref{coordinate}).
Recall from \cite{Pi1} that $k_1=k_1(\eta)$ and $k_2=k_2(\zeta)$ have triangular
factorizations of the following special forms:

\begin{equation}\label{k1formula}k_1=\left(\begin{matrix} 1&0\\
y^{*}&1\end{matrix} \right)\left(\begin{matrix} a_1&0\\
0&a_1^{-1}\end{matrix} \right)\left(\begin{matrix} \alpha_1&\beta_1\\
\gamma_1&\delta_1\end{matrix} \right),\end{equation}
and
\begin{equation}\label{k2formula}k_2=\left(\begin{matrix} 1&x^{*}\\
0&1\end{matrix} \right)\left(\begin{matrix} a_2&0\\
0&a_2^{-1}\end{matrix} \right)\left(\begin{matrix} \alpha_2&\beta_2\\
\gamma_2&\delta_2\end{matrix} \right)\end{equation}
where for example $y=\sum_{n\ge 0}y_nz^n$ and $x=\sum_{n\ge 1} x_nz^n$ are holomorphic functions
in $\Delta$, with appropriate boundary behavior, depending the smoothness of $g$.

As in \cite{Pi1}, given these triangular factorizations for $k_1$
and $k_2$, we can derive the triangular factorization for $g$ as
follows:
$$g=\left(\begin{matrix} \alpha_1&\beta_1\\
\gamma_1&\delta_1\end{matrix}
\right)^*\left(\begin{matrix}1&Y\\0&1\end{matrix}\right)\left(\begin{matrix}a_1a_2e^{\chi_{-}+\chi_0+\chi_{+}}&0
\\0&(a_1a_2e^{\chi_{-}+\chi_0+\chi_{+}})^{-1}\end{matrix}\right)
\left(\begin{matrix} 1&X^{*}\\
0&1\end{matrix} \right)\left(\begin{matrix} \alpha_2&\beta_2\\
\gamma_2&\delta_2\end{matrix} \right)$$

$$=\left(\begin{matrix} \alpha_1^*&\gamma_1^*\\
\beta_1^*&\delta_1^*\end{matrix}
\right)\left(\begin{matrix}e^{\chi_{-}}&0
\\0&e^{-\chi_-}\end{matrix}\right)$$
$$\left(\begin{matrix}1&e^{-2\chi_-}Y\\0&1\end{matrix}\right)
\left(\begin{matrix}a_1a_2e^{\chi_0}&0
\\0&(a_1a_2e^{\chi_0})^{-1}\end{matrix}\right)
\left(\begin{matrix} 1&e^{2\chi_{+}}X^{*}\\
0&1\end{matrix} \right)\left(\begin{matrix}e^{\chi_{+}}&0
\\0&e^{-\chi_{+}}\end{matrix}\right)\left(\begin{matrix} \alpha_2&\beta_2\\
\gamma_2&\delta_2\end{matrix} \right)$$
where $Y=a_1^{2}y$ and $X=a_2^{-2}x$.

The product of the middle
three factors is upper triangular, and it is easy to find its
triangular factorization.  Thus $g=l(g)m(g)a(g)u(g)$, where
\begin{equation}\label{lmatrix}l(g)=\left(\begin{matrix}l_{11}&l_{12}\\
l_{21}&l_{22}\end{matrix}\right)=\left(\begin{matrix} \alpha_1^*&\gamma_1^*\\
\beta_1^*&\delta_1^*\end{matrix}
\right)\left(\begin{matrix}e^{\chi_{-}}&0
\\0&e^{\chi^*}\end{matrix}\right)
\left(\begin{matrix}1&(e^{-2\chi_-}Y+(a_1a_2)^2e^{2(\chi_0+\chi_{+})}X^*)_-\\0&1\end{matrix}\right)
\end{equation}
\begin{equation}\label{ma}m(g)=\left(\begin{matrix}e^{\chi_0}&0\\0& e^{-\chi_0}\end{matrix}\right),
\quad a(g)=
\left(\begin{matrix}a_1a_2&0\\0&(a_1a_2)^{-1}\end{matrix}\right)\end{equation}
\begin{equation}\label{umatrix}u(g)=\left(\begin{matrix}u_{11}&u_{12}\\
u_{21}&u_{22}\end{matrix}\right)=\left(\begin{matrix} 1&((a_1a_2)^{-2}e^{-2(\chi_{-}+\chi_0)}Y+e^{2\chi_{+}}X^{*})_+\\
0&1\end{matrix} \right)\left(\begin{matrix}e^{\chi_{+}}&0
\\0&e^{-\chi_{+}}\end{matrix}\right)\left(\begin{matrix} \alpha_2&\beta_2\\
\gamma_2&\delta_2\end{matrix} \right)\end{equation}

In particular
\begin{equation}\label{lfactor1}\left(\begin{matrix}l_{11}\\
l_{21}\end{matrix}\right)=e^{\chi_{-}}\left(\begin{matrix}\alpha_1^*\\\beta_1^*\end{matrix}\right)\quad (a,b)=a_1(\alpha_1,\beta_1)$$
$$\left(u_{21},u_{22}\right)=
e^{-\chi_{+}}\left(\begin{matrix}\gamma_2,\delta_2\end{matrix}\right)\text{ and } (c,d)=a_2^{-1}(\gamma_2,\delta_2)\end{equation}
Therefore

\begin{equation}\label{lfactor2}\mathbb P\left(\begin{matrix}l_{11}^*\\
l_{21}^*\end{matrix}\right)=\mathbb
P\left(\begin{matrix}\alpha_1\\\beta_1\end{matrix}\right)=\mathbb
P\left(\begin{matrix}a\\b\end{matrix}\right)$$ and
$$
\mathbb P\left(u_{21},u_{22}\right)=
\mathbb
P\left(\gamma_2,\delta_2\right)=
\mathbb
P\left(c,d\right)\end{equation}

\section{Finding the Root Subgroup Factors: The $SU(2)$ Case}

Suppose that $g:S^1\to SU(2)$ (in what follows we will suppress mention of the necessary degree of smoothness
of $g$; but everything works if $g\in W^{1/2,L^2}(S^1,SU(2))$).

If $A(g)$ is invertible, then $g$ has a unique Birkhoff (or Riemann-Hilbert) factorization
factorization
$$g=g_- g_0 g_+$$
$g\in H^0(\Delta^*,G)$, $g_0\in G$, and $g_+\in H^0(\Delta,G)$, where
\begin{equation}\label{factorformula}(g_0g_+)^{-1}=[A(g)^{-1}\left(\begin{matrix}1\\0\end{matrix}\right),
A(g)^{-1}\left(\begin{matrix}0\\1\end{matrix}\right)]\end{equation}
where $A(g)$ is the (block) Toeplitz operator corresponding to $g$ (see Section 1 of \cite{Pi1}).

If $A(g)$ is invertible, and $g_0$ has a triangular
factorization in the finite dimensional sense, then $g$ has a
triangular factorization $g=lmau$. In this case we know that there also exists a root subgroup factorization
$$g=k_1(\eta)^* \left(\begin{matrix}e^{\chi}&0\\0&e^{-\chi}\end{matrix}\right) k_2(\zeta)$$
The key to solving for the $\eta$, $\chi$, and $\zeta$ factors is the following

\begin{theorem} (a) For $k_1(\eta)$,
$$b/a=\beta_1/\alpha_1=l_{21}^*/l_{11}^*$$
This meromorphic function in $\Delta$ has Taylor series $\sum_{n=0}^{\infty}\psi_nz^n$
where $\psi_n$ is the sum of terms
$$\psi_n= (-1)^r(-\overline \eta_{i_{0}})\left(\eta_{j_1}(-\overline \eta_{i_1})\right)...\left(\eta_{j_r}(-\overline \eta_{i_r})\right) $$
where $j_s<i_s$ and $j_s\le i_{s-1}$ for $s=1,..,r$, and $\sum_{s=1}^{r+1} i_s -\sum_{s=1}^r j_s=n$; in particular
$$\xi_n=(-\overline{\eta}_n)\prod_{s=1}^{n-1}(1+|\eta_s|^2)+polynomial(\eta_s,\overline{\eta}_s,s<n)$$

(b) For $k_2(\zeta)$,
$$c/d=\gamma_2/\delta_2=u_{21}/u_{22}=(g_+)_{21}/(g_+)_{22}$$
This meromorphic function in $\Delta$ has Taylor series $\sum_{n=1}^{\infty}\xi_nz^n$
where $\xi_n$ is the sum of terms
$$(-1)^r(-\overline \zeta_{i_{0}})\left(\zeta_{j_1}(-\overline \zeta_{i_1})\right)...\left(\zeta_{j_r}(-\overline \zeta_{i_r})\right) $$
where $j_s<i_s$ and $j_s\le i_{s-1}$ for $s=1,..,r$, and $\sum_{s=1}^{r+1} i_s -\sum_{s=1}^r j_s=n$; in particular
$$\xi_n=(-\overline{\zeta}_n)\prod_{s=1}^{n-1}(1+|\zeta_s|^2)+polynomial(\zeta_s,\overline{\zeta}_s,s<n)$$
\end{theorem}

For example
$$b/a=(-\overline{\eta}_0)+(-\overline{\eta}_1)(1+|\eta_0|^2)z+\left( (-\overline{\eta}_2)(1+|\eta_0|^2)(1+|\eta_1|^2)+...\right)z^2+... $$
and
$$c/d=(-\overline{\zeta}_1) z+(-\overline{\zeta}_2)(1+|\zeta_1|^2)z^2+\left((-\overline{\zeta}_3)(1+|\zeta_1|^2)(1+|\zeta_2|^2)+
(-\zeta_1\overline{\zeta}_2^2)(1+|\zeta_1|^2)\right)z^3$$
$$+((-\overline{\zeta}_4)(1+|\zeta_1|^2)(1+|\zeta_2|^2)(1+|\zeta_3|^2)
+(1+|\zeta_1|^2)(\zeta_2\overline{\zeta}_3^2(1+|\zeta_2|^2)$$
$$+2\zeta_1\overline{\zeta}_2\overline{\zeta}_3(1+|\zeta_2|^2)
+\overline{\zeta}_1^2\overline{\zeta}_2^3) )z^4+...$$ For this latter sum, if $n=2$, the terms in $\xi_2$ are $-\overline \zeta_2$ and $ (-1)(-\overline \zeta_1)\zeta_1(-\overline \zeta_2)$.

\begin{remarks} (a) Note that we can solve for the $\zeta$ variables using the coefficients of $(g_+)_{21}/(g_+)_{22}$.
We only need the Riemann-Hilbert factorization, not the full triangular factorization, in order to find $\zeta$.

(b) The finiteness of the formulas in the theorem contrasts sharply with the infinite formulas for the coefficients of
the terms $\gamma_2,\delta_2$ and so on; see Lemma \ref{keylemma} below.

(c) It would be nice to find a quicker route to finding $\chi$. It is essentially determined by either of the formulas
$$Re(\chi_+)=-log(a_1)-\frac12 log(\vert l_{11} \vert ^2+ \vert l_{21} \vert
^2))$$
$$=log(a_2)-\frac12 log(\vert u_{21} \vert ^2+\vert u_{22} \vert
^2) $$ where $a_1$ ($a_2$) is determined by $\eta$ ( $\zeta$, respectively); see \cite{Pi1}. Note that for the equality of
these two formulas, it is essential that the loop is unitary. For nonunitary loops these formulas diverge, and it is in
this sense that we said in the introduction, whereas $\eta$ and $\zeta$ extend naturally to the top stratum of the formal
completion of the complex loop group, $\chi$ does not have a preferred extension even to ordinary complex loops; $\chi$
is much more complicated from an analytic point of view.

\end{remarks}

\begin{proof} We will prove part (b); part (a) is proven in the same way.

We first consider the equality of the meromorphic functions, $b/a$, $\gamma_2/\delta_2$, $u_{21}/u_{22}$, and $(g_+)_{21}/(g_+)_{22}$.
The first two equalities follow from (\ref{lfactor1}).  As we noted before the statement of the theorem, a
triangular factorization implies a Riemann-Hilbert factorization:
If $g=lmau$, then $g=g_-g_0g_+$, where
$$g_-=\left(\begin{matrix}l_{11}-l_{12}l_{21}(\infty)&l_{12}\\
l_{21}-l_{22}l_{21}(\infty)&l_{22}\end{matrix}\right)$$
$$g_0=\left(\begin{matrix}1&0\\
l_{21}(\infty)&1\end{matrix}\right)\left(\begin{matrix}a_1a_2e^{\chi_0}&0\\
0&(a_1a_2e^{\chi_0})^{-1}\end{matrix}\right)\left(\begin{matrix}1&u_{12}(0)\\
0&1\end{matrix}\right)$$
$$g_+=\left(\begin{matrix}u_{11}-u_{21}u_{12}(0)&u_{12}-u_{22}u_{12}(0)\\
u_{21}&u_{22}\end{matrix}\right)$$
In particular
$(g_+)_{21}=u_{21}$ and $(g_+)_{22}=u_{22}$. This obviously implies the third equality.

We now turn to proving the formula for $\xi_n$. Using $\xi=\gamma_2/\delta_2$, we need to show that for the product\begin{equation}\label{2product2}\left(\begin{matrix}\delta_2^*(z)&-\gamma_2^*(z)\\
\gamma_2(z)&\delta_2(z)\end{matrix}\right)=\lim_{n\to\infty}\mathbf a(\zeta_n)\left(\begin{matrix} 1&\zeta_nz^{-n}\\
-\overline{\zeta}_nz^n&1\end{matrix} \right)..\mathbf
a(\zeta_1)\left(\begin{matrix} 1&
\zeta_1z^{-1}\\
-\overline{\zeta}_1z&1\end{matrix} \right),\end{equation}
$\gamma_2/\delta_2$ has the claimed form. Our strategy is completely straightforward: we will first recall
the formulas for the coefficients of $\gamma_2$ and $\delta_2$, and we will then calculate the Taylor series for the quotient.

The following is from \cite{Pi1}, which is reminiscent of the Pauli exclusion principle:

\begin{lemma}\label{keylemma} For the product
\begin{equation}\label{2product}=\left(\begin{matrix}\delta_2^*(z)&-\gamma_2^*(z)\\
\gamma_2(z)&\delta_2(z)\end{matrix}\right)=\lim_{n\to\infty}\mathbf a(\zeta_n)\left(\begin{matrix} 1&\zeta_nz^{-n}\\
-\overline{\zeta}_nz^n&1\end{matrix} \right)..\mathbf
a(\zeta_1)\left(\begin{matrix} 1&
\zeta_1z^{-1}\\
-\overline{\zeta}_1z&1\end{matrix} \right),\end{equation}
where
$$\gamma_2(z)=\sum_{n=1}^{\infty}\gamma_{2,n}z^n,$$
$$\gamma_{2,n}=\sum (-\overline{\zeta}_{i_1})\zeta_{j_1}...(-\overline{\zeta}_{
i_r})\zeta_{j_r}(-\overline{\zeta}_{i_{r+1}}),$$ the sum over
multiindices satisfying
$$0<i_1<j_1<..<j_r<i_{r+1},\quad\sum i_{*}-\sum j_{*}=n,$$
and
$$\delta_2(z)=1+\sum_{n=1}^{\infty}\delta_{2,n}z^n,$$
$$\delta_{2,n}=\sum\zeta_{j_1}(-\overline{\zeta}_{i_1})...\zeta_{j_r}(
-\overline{\zeta}_{i_r}),$$ the sum over multiindices satisfying
$$0<j_1<i_1<..<i_r,\quad\sum (i_{*}-j_{*})=n$$
\end{lemma}

To simplify notation, let $\gamma:=\gamma_2$, and write $\delta_2:=1+\delta$.
Then $\gamma_2/\delta_2=\gamma(1-\delta+\delta^2-..$, and the nth coefficient of $\gamma_2/\delta_2$ equals
$$\gamma_n-(\gamma\delta)_n+(\gamma\delta^2)_n-..+(-1)^{n-1}(\gamma\delta^{n-1})_n$$
Each of the terms in this sum, according to the Lemma, has an expression as an infinite sum. According to the statement
of the theorem, all but finitely many of these terms cancel out.

To explain in a leisurely way how this comes about, first consider
\begin{equation}\label{sumofterms}(\gamma\delta)_n=\sum_{k=1}^{n-1}\gamma_{k}\delta_{n-k}\end{equation}
This is a sum of terms of the form
$$\left(-\overline{\zeta}_{i_1})\zeta_{j_1}...(-\overline{\zeta}_{
i_r})\zeta_{j_r}(-\overline{\zeta}_{i_{r+1}})\right)\left(\zeta_{j'_1}(-\overline{\zeta}_{i'_1})...\zeta_{j'_{r'}}(
-\overline{\zeta}_{i'_{r'}})\right)$$
where
$$0<i_1<j_1<..<j_r<i_{r+1}, \sum i_{*}-\sum j_{*}=k, 0<j'_1<i'_1<..<i'_{r'},\sum i'_{*}-\sum j'_{*}=n-k$$
If $i_{r+1}<j'_1$, then this product will exactly cancel with a term in the corresponding sum for $\gamma_n$; it is in some sense
obeying a Pauli exclusion principle.  The only term in the sum for $\gamma_n$ that is not canceled is $-\overline{\zeta}_n$; this is the only term which cannot be broken into two terms as in the sum for $\gamma\delta$. If $j'_1\le i_{r+1}$, then we keep this term; however,
we will see than many of these terms are canceled by subsequent terms appearing in the sum (\ref{sumofterms}).

Consider $(\gamma\delta^s)_n$. This is a sum of terms of the form
\begin{equation}\label{term1}\left(-\overline{\zeta}_{i_1})\zeta_{j_1}...(-\overline{\zeta}_{
i_r})\zeta_{j_r}(-\overline{\zeta}_{i_{r+1}}\right) \times
\left(\zeta_{j^1_1}(-\overline{\zeta}_{i^1_1})...\zeta_{j^1_{r^1}}(
-\overline{\zeta}_{i^1_{r^1}})\right)\times ... \times \left(\zeta_{j^s_1}(-\overline{\zeta}_{i^s_1})...\zeta_{j^s_{r^s}}(
-\overline{\zeta}_{i^s_{r^s}})\right)\end{equation}
where each of the factors separated by $\times$ (the first factor comes
from $\gamma$, and the other $s$ factors come from $\delta^s$) satisfy the appropriate constraints in the Lemma (in particular the sum
of the $i$ indices minus the sum of the $j$ indices equals $n$).

At one extreme, it may happen that all of the indices in (\ref{term1}) are increasing, i.e. $i_{r+1}<j^1_1$ and $i^{s'}_{r^{s'}}<j^{s'+1}_1$ for $s'=1,..,s-1$. By removing some of the $\times$, we see that this product will have occurred in
all of the preceding terms $(\gamma\delta^{s'})_n$, $s'=0,..,s-1$ (which occur with alternating signs). Similarly if
$r>0$ or $r^j>1$ for some $j$, then we can insert $\times$ and this more finely factored product will occur in some subsequent terms $(\gamma\delta^{s'})_n$, for $s<s'$; the largest such $s'$ is $S=r+r^1+..+r^s$, in which for each factor
of $\delta$, the corresponding factor
$$\zeta_{j_1}(-\overline{\zeta}_{i_1})...\zeta_{j_{r'}}(
-\overline{\zeta}_{i_{r'}})$$ is irreducible in the sense that it cannot be split into a product of two similar factors, i.e. $r'=1$.
The number of times the product (\ref{term1}) occurs in one of the terms in (\ref{sumofterms}) depends on the number of ways we can insert $\times$. Taking into account
the signs that occur in (\ref{sumofterms}), the coefficient of this product in (\ref{sumofterms}) is
$$\sum_{j=0}^S(-1)^j\left(\begin{matrix}S\\j\end{matrix}\right)=0$$ Thus this
product completely cancels out in the sum (\ref{sumofterms}).

At the opposite extreme, the term (\ref{term1}) may have the same form as in the statement of the theorem, i.e. $r=0$, and
$r^j=1$, $j=1,..,s$. In this case this term occurs in exactly one of the terms $(\gamma\delta^{s'})_n$. There is no cancelation.

In between these two extremes, we are considering a product for which, in the expression (\ref{term1}), there is a positive number of instances when $j^{s''}_{1}\le i^{s''-1}_{r^{s''-1}}$ for some $s''=1,..,s-1$. As in the first extreme case, we can possibly remove some of the $\times$ to
see that this term occurs in earlier terms $(\gamma\delta^{s'})_n$ ($s'<s$), and we can possibly insert some $\times$ to see that it occurs in some
later terms $(\gamma\delta^{s'})_n$ ($s<s'$). For definiteness we can suppose that $s$ is as large as possible, i.e. that $r=0$
and each $r^j=1$, so that it is just a question of removing $\times$. If $s_0$ is the smallest $s'$ such that the product occurs
in $\gamma\delta^{s'}$, then $s_0<s$ (because we are not in the second extreme case) and the coefficient of this product in (\ref{sumofterms}) is
$$\sum_{j=s_0}^s(-1)^j\left(\begin{matrix}s-s_0\\j\end{matrix}\right)=0$$
Thus this product completely cancels out.

If we multiply out $(-\overline \zeta_n)\prod(1+\zeta_i\overline \zeta_i)$, then we see that each of the terms does occur in
the sum in part (b). This proves the last claim in part (b).

\end{proof}

\subsection{Solving for $\eta$ and $\zeta$}

\begin{corollary} (a) $\eta_i$ is a rational function of the coefficients $\psi_{i'}$, and in turn the coefficients $\psi$
are polynomials in the coefficients of $l_{21}$ and $l_{22}$.

(b) $\zeta_k$ is a rational function of the coefficients $\xi_{k'}$, and in turn the coefficients $\xi$
are polynomials in the coefficients of $u_{21}$ and $u_{22}$.

\end{corollary}

Unfortunately (based on Maple calculations) it appears hopeless to find a closed form expression for the $\zeta$ variables in terms of the $\xi$ variables.

\section{$\eta$ and $\zeta$ as Functions on the Formal Completion}

The formal completion of the loop group $L\dot G$ is defined by
$$\mathbf LG=G(\mathbb C((z^{-1})))\times_{G(\mathbb C[z,z^{-1}])}G(
\mathbb C((z))),$$ where $\mathbb C((z))$ is the field of formal
Laurent series $\sum a_nz^n$, $a_n=0$ for $n<<0$. There is a generalized
Birkhoff decomposition
$$\mathbf LG=\bigsqcup_{\lambda\in Hom(S^1,T)}\Sigma_{\lambda}^{\mathbf LG},\quad
\Sigma_{\lambda}^{\mathbf LG}=G(\mathbb
C[[z^{-1}]])\cdot\lambda\cdot G(\mathbb C[[z]])$$ where $\mathbb
C[[z ]]$ denotes formal power series in $z$. This decomposition
reduces to the usual Birkhoff decomposition for the smooth loop group $L\dot G$.

For many purposes one is primarily interested in the top stratum corresponding to $\lambda=1$ (which is open and dense).
For $g$ in the top stratum, there is a unique formal Riemann-Hilbert factorization
$$g=g_{-}\cdot g_0\cdot g_{+}$$
where $g_-\in G(\mathbb
C[[z^{-1}]])$, $g_-(\infty)=1$ $g_0 \in G$, and $g_+\in G(\mathbb C[[z]])$, $g_+(0)=1$.
There are bijective correspondences
$$\{g_-\in G(\mathbb
C[[z^{-1}]]): g_-(\infty)=1\}\leftrightarrow \mathfrak g(\mathbb
C[[z^{-1}]]): g_- \leftrightarrow \theta_-=(\partial g_-)g_-^{-1}$$ and
$$\{g_+\in G(\mathbb
C[[z]]): g_+(0)=1\}\leftrightarrow \mathfrak g(\mathbb
C[[z]]): g_+ \leftrightarrow \theta_+=g_+^{-1}(\partial g_+$$
The factors $g_+,g_0,g_-$ ($\theta_-,g_0,\theta_+$) are referred to as Riemann-Hilbert coordinates
(linear Riemann-Hilbert
coordinates, respectively) for $g\in \mathbb L \dot G$.

If $g_0$ additionally has a triangular factorization, then $g$
has a unique formal triangular factorization
\begin{equation}\label{trifactorization2}g=l\cdot m \cdot a\cdot u \end{equation}
where $m\in T$ (the diagonal
torus in $SU(2)$),
$a\in A=exp(\mathbb R\left(\begin{matrix} 1&0\\
0&-1\end{matrix} \right))$, $l\in \mathcal N^{-}$, the (profinite
nilpotent) group consisting of formal power series in $z^{-1}$,
$$l=\left(\begin{matrix} 1+\sum_{j=1}^{\infty}A_jz^{-j}&\sum_{j=1}^{\infty}B_jz^{
-j}\\
\sum_{j=0}^{\infty}C_jz^{-j}&1+\sum_{j=1}^{\infty}D_jz^{-j}&\end{matrix}
\right ),$$ with $det(l)=1$ (as a formal power series in
$z^{-1}$), and $u\in \mathcal N^{+}$, the group consisting of
formal power series in $z$,
$$u=\left(\begin{matrix} 1+\sum_{j=1}^{\infty}a_jz^j&\sum_{j=0}^{\infty}b_jz^j\\
\sum_{j=1}^{\infty}c_jz^j&1+\sum_{j=1}^{\infty}d_jz^j&\end{matrix}
\right),$$ with $det(u)=1$.

Exactly as in the previous subsection, we can say that $\eta$ and $\zeta$ extend in a natural way
to functions on the set of $g\in \mathbf L G$ having a triangular factorization.

\begin{corollary} In reference to the formal completion,
(a) $\eta_i$ is a rational function of the coefficients of $l_{21}$ and $l_{22}$, which are in
turn polynomials in the coefficients of $\theta_-$.

(b) $\zeta_k$ is a rational function of the coefficients of $u_{21}$ and $u_{22}$, which are in turn polynomials in the coefficients
of $\theta_+$.

\end{corollary}

This has important measure-theoretic implications, which we will pursue elsewhere.

\end{document}